\documentclass[12pt]{article}

\catcode`\@=11

\thicklines
\newskip\Einheit \Einheit=.6cm
\newcount\xcoord \newcount\ycoord
\newdimen\xdim \newdimen\ydim \newdimen\PfadD@cke \newdimen\Pfadd@cke
\def\PfadDicke#1{\PfadD@cke#1 \divide\PfadD@cke by2 
\Pfadd@cke\PfadD@cke \multiply\PfadD@cke by2}
\long\def\LOOP#1\REPEAT{\def\BODY{#1}\ITERATE}
\def\ITERATE{\BODY \let\next\ITERATE \else\let\next\relax\fi \next}
\let\REPEAT=\fi
\def\Punkt{\hbox{\raise-2pt\hbox to0pt{\hss\scriptsize$\bullet$\hss}}}

\def\DuennPunkt(#1,#2){\unskip
  \raise#2 \Einheit\hbox to0pt{\hskip#1 \Einheit
          \raise-1.5pt\hbox to0pt{\hss\tiny$\bullet$\hss}\hss}}
		  
\def\NormalPunkt(#1,#2){\unskip
  \raise#2 \Einheit\hbox to0pt{\hskip#1 \Einheit
          \raise-3pt\hbox to0pt{\hss\large$\bullet$\hss}\hss}}
\def\DickPunkt(#1,#2){\unskip
  \raise#2 \Einheit\hbox to0pt{\hskip#1 \Einheit
          \raise-4pt\hbox to0pt{\hss\Large$\bullet$\hss}\hss}}
\def\Kreis(#1,#2){\unskip
  \raise#2 \Einheit\hbox to0pt{\hskip#1 \Einheit
          \raise-4pt\hbox to0pt{\hss\Large$\circ$\hss}\hss}}
\def\Diagonale(#1,#2)#3{\unskip\leavevmode
  \xcoord#1\relax \ycoord#2\relax
      \raise\ycoord \Einheit\hbox to0pt{\hskip\xcoord \Einheit
         \unitlength\Einheit
         \line(1,1){#3}\hss}}
\def\AntiDiagonale(#1,#2)#3{\unskip\leavevmode
  \xcoord#1\relax \ycoord#2\relax \advance\xcoord by -0.05\relax
      \raise\ycoord \Einheit\hbox to0pt{\hskip\xcoord \Einheit
         \unitlength\Einheit
         \line(1,-1){#3}\hss}}
\def\Pfad(#1,#2),#3\endPfad{\unskip\leavevmode
  \xcoord#1 \ycoord#2 \thicklines\ZeichnePfad#3\endPfad\thinlines}
\def\ZeichnePfad#1{\ifx#1\endPfad\let\next\relax
  \else\let\next\ZeichnePfad
    \ifnum#1=1
      \raise\ycoord \Einheit\hbox to0pt{\hskip\xcoord \Einheit
         \vrule height\Pfadd@cke width1 \Einheit depth\Pfadd@cke\hss}%
      \advance\xcoord by 1
    \else\ifnum#1=2
      \raise\ycoord \Einheit\hbox to0pt{\hskip\xcoord \Einheit
        \hbox{\hskip-\PfadD@cke\vrule height1 \Einheit 
width\PfadD@cke depth0pt}\hss}%
      \advance\ycoord by 1
    \else\ifnum#1=3
      \raise\ycoord \Einheit\hbox to0pt{\hskip\xcoord \Einheit
         \unitlength\Einheit
         \line(1,1){1}\hss}
      \advance\xcoord by 1
      \advance\ycoord by 1
    \else\ifnum#1=4
      \raise\ycoord \Einheit\hbox to0pt{\hskip\xcoord \Einheit
         \unitlength\Einheit
         \line(1,-1){1}\hss}
      \advance\xcoord by 1
      \advance\ycoord by -1
    \fi\fi\fi\fi
  \fi\next}
\def\hSSchritt{\leavevmode\raise-.4pt\hbox 
to0pt{\hss.\hss}\hskip.2\Einheit
  \raise-.4pt\hbox to0pt{\hss.\hss}\hskip.2\Einheit
  \raise-.4pt\hbox to0pt{\hss.\hss}\hskip.2\Einheit
  \raise-.4pt\hbox to0pt{\hss.\hss}\hskip.2\Einheit
  \raise-.4pt\hbox to0pt{\hss.\hss}\hskip.2\Einheit}
\def\vSSchritt{\vbox{\baselineskip.2\Einheit\lineskiplimit0pt
\hbox{.}\hbox{.}\hbox{.}\hbox{.}\hbox{.}}}
\def\DSSchritt{\leavevmode\raise-.4pt\hbox to0pt{%
  \hbox to0pt{\hss.\hss}\hskip.2\Einheit
  \raise.2\Einheit\hbox to0pt{\hss.\hss}\hskip.2\Einheit
  \raise.4\Einheit\hbox to0pt{\hss.\hss}\hskip.2\Einheit
  \raise.6\Einheit\hbox to0pt{\hss.\hss}\hskip.2\Einheit
  \raise.8\Einheit\hbox to0pt{\hss.\hss}\hss}}
\def\dSSchritt{\leavevmode\raise-.4pt\hbox to0pt{%
  \hbox to0pt{\hss.\hss}\hskip.2\Einheit
  \raise-.2\Einheit\hbox to0pt{\hss.\hss}\hskip.2\Einheit
  \raise-.4\Einheit\hbox to0pt{\hss.\hss}\hskip.2\Einheit
  \raise-.6\Einheit\hbox to0pt{\hss.\hss}\hskip.2\Einheit
  \raise-.8\Einheit\hbox to0pt{\hss.\hss}\hss}}
\def\SPfad(#1,#2),#3\endSPfad{\unskip\leavevmode
  \xcoord#1 \ycoord#2 \ZeichneSPfad#3\endSPfad}
\def\ZeichneSPfad#1{\ifx#1\endSPfad\let\next\relax
  \else\let\next\ZeichneSPfad
    \ifnum#1=1
      \raise\ycoord \Einheit\hbox to0pt{\hskip\xcoord \Einheit
         \hSSchritt\hss}%
      \advance\xcoord by 1
    \else\ifnum#1=2
      \raise\ycoord \Einheit\hbox to0pt{\hskip\xcoord \Einheit
        \hbox{\hskip-2pt \vSSchritt}\hss}%
      \advance\ycoord by 1
    \else\ifnum#1=3
      \raise\ycoord \Einheit\hbox to0pt{\hskip\xcoord \Einheit
         \DSSchritt\hss}
      \advance\xcoord by 1
      \advance\ycoord by 1
    \else\ifnum#1=4
      \raise\ycoord \Einheit\hbox to0pt{\hskip\xcoord \Einheit
         \dSSchritt\hss}
      \advance\xcoord by 1
      \advance\ycoord by -1
    \fi\fi\fi\fi
  \fi\next}
\def\Koordinatenachsen(#1,#2){\unskip
 \hbox to0pt{\hskip-.5pt\vrule height#2 \Einheit width.5pt depth1 
\Einheit}%
 \hbox to0pt{\hskip-1 \Einheit \xcoord#1 \advance\xcoord by1
    \vrule height0.25pt width\xcoord \Einheit depth0.25pt\hss}}
\def\Koordinatenachsen(#1,#2)(#3,#4){\unskip
 \hbox to0pt{\hskip-.5pt \ycoord-#4 \advance\ycoord by1
    \vrule height#2 \Einheit width.5pt depth\ycoord \Einheit}%
 \hbox to0pt{\hskip-1 \Einheit \hskip#3\Einheit 
    \xcoord#1 \advance\xcoord by1 \advance\xcoord by-#3 
    \vrule height0.25pt width\xcoord \Einheit depth0.25pt\hss}}
\def\Gitter(#1,#2){\unskip \xcoord0 \ycoord0 \leavevmode
  \LOOP\ifnum\ycoord<#2
    \loop\ifnum\xcoord<#1
      \raise\ycoord \Einheit\hbox to0pt{\hskip\xcoord 
\Einheit\Punkt\hss}%
      \advance\xcoord by1
    \repeat
    \xcoord0
    \advance\ycoord by1
  \REPEAT}
\def\Gitter(#1,#2)(#3,#4){\unskip \xcoord#3 \ycoord#4 \leavevmode
  \LOOP\ifnum\ycoord<#2
    \loop\ifnum\xcoord<#1
      \raise\ycoord \Einheit\hbox to0pt{\hskip\xcoord 
\Einheit\Punkt\hss}%
      \advance\xcoord by1
    \repeat
    \xcoord#3
    \advance\ycoord by1
  \REPEAT}
\def\Label#1#2(#3,#4){\unskip \xdim#3 \Einheit \ydim#4 \Einheit
  \def\lo{\advance\xdim by-.5 \Einheit \advance\ydim by.5 \Einheit}%
  \def\llo{\advance\xdim by-.25cm \advance\ydim by.5 \Einheit}%
  \def\loo{\advance\xdim by-.5 \Einheit \advance\ydim by.25cm}%
  \def\o{\advance\ydim by.25cm}%
  \def\ro{\advance\xdim by.5 \Einheit \advance\ydim by.5 \Einheit}%
  \def\rro{\advance\xdim by.25cm \advance\ydim by.5 \Einheit}%
  \def\roo{\advance\xdim by.5 \Einheit \advance\ydim by.25cm}%
  \def\l{\advance\xdim by-.30cm}%
  \def\r{\advance\xdim by.30cm}%
  \def\lu{\advance\xdim by-.5 \Einheit \advance\ydim by-.6 \Einheit}%
  \def\llu{\advance\xdim by-.25cm \advance\ydim by-.6 \Einheit}%
  \def\luu{\advance\xdim by-.5 \Einheit \advance\ydim by-.30cm}%
  \def\u{\advance\ydim by-.30cm}%
  \def\ru{\advance\xdim by.5 \Einheit \advance\ydim by-.6 \Einheit}%
  \def\rru{\advance\xdim by.25cm \advance\ydim by-.6 \Einheit}%
  \def\ruu{\advance\xdim by.5 \Einheit \advance\ydim by-.30cm}%
  #1\raise\ydim\hbox to0pt{\hskip\xdim
     \vbox to0pt{\vss\hbox to0pt{\hss$#2$\hss}\vss}\hss}%
}
\catcode`\@=12

\usepackage{pstricks,pst-node,pst-tree}
\usepackage{amsmath,amsthm}
\usepackage{color}

\usepackage{xspace}
\usepackage[colorlinks=true,
linkcolor=green,
filecolor=brown,
citecolor=green]{hyperref}

\def\v{\vert}
\def\ve{\varepsilon}

\def\p{\pstree}
\def\m{\ensuremath{\mathcal M}\xspace}
\def\scr{\small}

\begin{document}
\newtheorem{bij}{Bijection}
\newtheorem{theorem}{Theorem}
\newtheorem{prop}{Proposition}
\newtheorem{cor}{Corollary}
\begin{center}
{\Large
Some Bijections for Restricted Motzkin Paths                         \\ 
}
\vspace{5mm}
DAVID CALLAN  \\
Department of Statistics  \\
University of Wisconsin-Madison  \\
1210 W. Dayton St   \\
Madison, WI \ 53706-1693  \\
{\bf callan@stat.wisc.edu}  \\
\vspace{2mm}
July 19 2004
\end{center}

\vspace{5mm}

%%%%%%%%%%%%%%%%%%%%%%%%%%%%%%%%%%%%%%%%%%%%%%%%%%%%%%%%%%%%%%%%%%%%%%%%%%%%%%%%%%%%%%%%%%%%%%%
%                                                                                             %
%   \centerline{\textbf{Abstract}}                                                            %
%   We give several bijections among restricted Motzkin paths, explaining why                 %
%   various parameters on these paths are equidistributed.  For example, the                  %
%   number of doublerise-free Motzkin paths of length $n$ is the same as the number of        %
%   peak-free Motzkin paths of length $n+1$ and the parameter ``number of doublefalls'' has   %
%   the same distribution on the former set as ``number of valleys'' does on the              %
%   latter.  The bijections are most easily presented recursively but we also                 %
%   give explicit descriptions using the notion of Motzkin tree.                              %
%                                                                                             %
%%%%%%%%%%%%%%%%%%%%%%%%%%%%%%%%%%%%%%%%%%%%%%%%%%%%%%%%%%%%%%%%%%%%%%%%%%%%%%%%%%%%%%%%%%%%%%%

There is a classic correspondence between full binary trees on $2n$ 
edges and Dyck paths of $2n$ steps \cite[Ex. 6.19d, 6.19i]{ec2}: traverse the tree in preorder 
(counterclockwise from the root) and, as each edge is encountered for 
the first time, record an upstep for a left edge and a downstep for a 
right edge. This correspondence readily extends to Motzkin paths 
because a Motzkin path can be viewed as a Dyck path with a nonnegative 
label on each vertex recording the number of flatsteps at that 
location, and Dyck path vertices in left-to-right order correspond to 
tree vertices in preorder. So just label each vertex in the tree with 
the number of flatsteps at the corresponding location in the path as 
illustrated below.

\hspace*{10mm}\p[levelsep=30pt,nodesep=3pt]{\TR{\scr{3}}} {\p{\TR{\scr{1}}} 
{\TR{\scr{ 0}} \p{\TR{\scr{0 }}} {\p{\TR{\scr{0}}} {\TR{\scr{2}} \TR{\scr{0}}}\p{\TR{\scr{1}}} {\TR{\scr{0}} \TR{\scr{0}}}}}
 \TR{\scr{0}}}

% 1=step in x-direction, 2=step in y-direction, 3=upward diagonal step, 
% 4=downward diagonal step
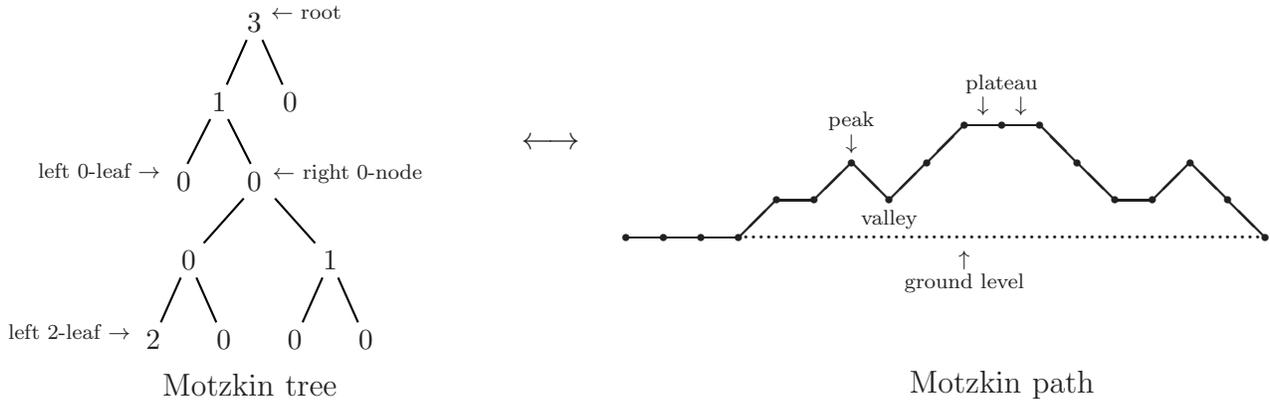
\vspace*{-48.5mm}

\Einheit=0.5cm
\[
\Label\o{\textrm{\scriptsize{left 0-leaf $\rightarrow$}}}(-14.0,1.3)
\Label\o{\textrm{\scriptsize{left 2-leaf $\rightarrow$}}}(-14.8,-3)
\Label\o{\textrm{\scriptsize{$\leftarrow$ right 0-node }}}(-7.3,1.2)
\Label\o{\textrm{\scriptsize{$\leftarrow$ root}}}(-8.5,5.5)
\Pfad(0,0),11131343311441344\endPfad
\SPfad(3,0),11111111111111\endSPfad
\DuennPunkt(0,0)
\DuennPunkt(1,0)
\DuennPunkt(2,0)
\DuennPunkt(3,0)
\DuennPunkt(4,1)
\DuennPunkt(5,1)
\DuennPunkt(6,2)
\DuennPunkt(7,1)
\DuennPunkt(8,2)
\DuennPunkt(9,3)
\DuennPunkt(10,3)
\DuennPunkt(11,3)
\DuennPunkt(12,2)
\DuennPunkt(13,1)
\DuennPunkt(14,1) 
\DuennPunkt(15,2)
\DuennPunkt(16,1)
\DuennPunkt(17,0) 
\Label\o{\longleftrightarrow}(-2,2)
\Label\o{\textrm{\scriptsize{peak}}}(6,2.6)
\Label\o{\textrm{\scriptsize{$\downarrow$}}}(6,2)
\Label\o{\textrm{\scriptsize{valley}}}(7,0)
\Label\o{\textrm{\scriptsize{plateau}}}(10,3.6)
\Label\o{\textrm{\scriptsize{$\downarrow$}}}(9.5,3)
\Label\o{\textrm{\scriptsize{$\downarrow$}}}(10.5,3)
\Label\u{\textrm{\scriptsize{ground level}}}(9,-.6)
\Label\u{\textrm{\scriptsize{$\uparrow$}}}(9,0)
%\Label\u{\longleftrightarrow}(-1,-3)
\Label\u{\textrm{Motzkin path}}(10,-3.3)
\Label\u{\textrm{Motzkin tree}}(-10,-3.3)
\]

\begin{center}
   \textbf{ Motzkin tree-path correspondence}
\end{center}

%\vspace*{-1mm}

The weight of the labeled tree is \#\,edges + sum of labels, and the 
length of a path is its number of steps. So weight 
of tree $\leftrightarrow$ length of path. We will call such a labeled 
tree of weight $n$ a Motzkin $n$-tree and a Motzkin path of length $n$  
a Motzkin $n$-path so that Motzkin $n$-trees correspond 
to Motzkin $n$-paths. (Motzkin $n$-trees are closely related to the 
$\{0,1,2\}$-trees of \cite{deutschMotzkin}.) We mostly follow the 
notation in \cite{deutschMotzkin}, in particular distinguishing 
between a node (non-root interior vertex) and a leaf. Thus the vertices 
of a Motzkin tree are partitioned into a root, a set of nodes, and a set 
of leaves. Each node and 
leaf is left or right according as it is a left or right child of its 
parent. A $k$-node is one whose label is $k$ and a positive node is 
one whose label is $\ge 1$. Similarly for leaves. The trivial tree has 
no leaves. In a nontrivial tree, the first and last 
leaf are as encountered in preorder, so the first leaf is left and 
the last is right. The level of a vertex is the 
length of the unique path joining it to the root. Every non-root 
vertex has a unique sibling---the other child of its parent. We use $U$ for 
upstep, $F$ for flatstep, and $D$ for downstep. A plateau in a path 
is a run of flatsteps that is either the entire path or of length $\ge 
1$ and preceded by an upstep and 
followed by a downstep.

We recall some obvious correspondences.
\begin{center}
\begin{tabular}{rl}
   \textbf{Motzkin tree} &  \textbf{Motzkin path}     \\
    %\hline
    root label & \#\,initial $F$s  \\
    %\hline
    first leaf & first peak or plateau  \\
    %\hline
    last leaf & terminal vertex  \\
    %\hline
    left $0$-node & doublerise ($UU$)  \\
    %\hline
    right  $0$-node & valley ($DU$)   \\
    %\hline
    left $0$-leaf & peak ($UD$)  \\
    %\hline
    \vspace*{-2mm}
    right $0$-leaf & doublefall ($DD$)  \\ 
   (except last leaf) &    \\
    %\hline 
    \vspace*{-2mm}
    level of first leaf & height of first peak  \\
      &  or plateau  \\
    %\hline 
    \vspace*{-2mm}
   level of last leaf & \#\,$D$s that return path \\
     &  to ground level \\
    %\hline
    \vspace*{-2mm}
    positive labels & plateau lengths  \\
   on left leaves &    \\
    %\hline
\end{tabular}
\end{center}

We use $\m_{n}(UD,DU)$ to denote the set of Motzkin $n$-paths that 
contain neither peaks nor valleys, and so on. Thus, for example $\m_{n}(UU)$
corresponds to the set of Motzkin $n$-trees in which each left node is 
positive.
 
Each of the following 5 bijections has both a recursive and an explicit 
description. 
The recursive specification depends on the first few 
steps and the first return 
to ground level; $\ve$ denotes the empty path, $R,S,T$ denote Motzkin 
paths. Motzkin trees facilitate the explicit description. 
The equivalence of the recursive and explicit descriptions can 
be proved by induction. Emeric Deutsch \cite{deutschUnpublished} 
found the recursive form of most of them.

\begin{bij}
    $\phi: \m_{n} \rightarrow \m_{n}$.
\end{bij}
Recursive:
\begin{eqnarray*}
    \phi(\ve) & = & \ve  \\
    \phi(FR) & = & F\phi(R)  \\    
    \phi(URDS) & = &   U\phi(S)D \phi(R)  
\end{eqnarray*}

Explicit: 
\[
\textrm{path\ } 
\xrightarrow[\textrm{above}]{\textrm{as}} 
\textrm{\ tree\ } 
\xrightarrow[\textrm{vertical}]{\textrm{flip in}} 
\textrm{\ tree\ } 
\xrightarrow[\textrm{above}]{\textrm{as}} 
\textrm{\ path}
\]

Example: (0 labels omitted)
% 1=step in x-direction, 2=step in y-direction, 3=upward diagonal step, 
% 4=downward diagonal step
\vspace*{2mm}
\Einheit=0.5cm
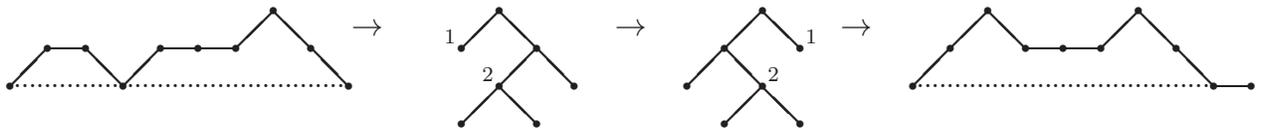
\[
\Pfad(-16,1),314311344\endPfad
\SPfad(-16,1),111111111\endSPfad
\Pfad(-4,2),344\endPfad
\Pfad(-4,0),33\endPfad
\Pfad(-3,1),4\endPfad
\Pfad(2,1),334\endPfad
\Pfad(3,2),4\endPfad
\Pfad(3,0),34\endPfad
\Pfad(8,1),334113441\endPfad
\SPfad(8,1),11111111\endSPfad
\DuennPunkt(-16,1)
\DuennPunkt(-15,2)
\DuennPunkt(-14,2)
\DuennPunkt(-13,1)
\DuennPunkt(-12,2)
\DuennPunkt(-11,2)
\DuennPunkt(-10,2)
\DuennPunkt(-9,3)
\DuennPunkt(-8,2)
\DuennPunkt(-7,1)
\DuennPunkt(-4,2)
\DuennPunkt(-4,0)
\DuennPunkt(-3,1)
\DuennPunkt(-3,3)
\DuennPunkt(-2,2)
\DuennPunkt(-2,0)
\DuennPunkt(-1,1)
\DuennPunkt(2,1)
\DuennPunkt(3,0)
\DuennPunkt(3,2)
\DuennPunkt(4,1)
\DuennPunkt(4,3)
\DuennPunkt(5,2)
\DuennPunkt(5,0) 
\DuennPunkt(8,1)
\DuennPunkt(9,2)
\DuennPunkt(10,3)
\DuennPunkt(11,2)
\DuennPunkt(12,2)
\DuennPunkt(13,2)
\DuennPunkt(14,3) 
\DuennPunkt(15,2)
\DuennPunkt(16,1)
\DuennPunkt(17,1) 
\Label\o{\rightarrow}(-6.5,2)
\Label\o{\rightarrow}(0.5,2)
\Label\o{\rightarrow}(6.5,2)
\Label\o{\textrm{{\scriptsize 2}}}(-3.3,0.8)
\Label\o{\textrm{{\scriptsize 1}}}(-4.3,1.8)
\Label\o{\textrm{{\scriptsize 2}}}(4.3,0.8)
\Label\o{\textrm{{\scriptsize 1}}}(5.3,1.8)
\] 

Consequence:\quad Since left $0$-nodes ($UU$) and right $0$-nodes ($DU$) 
are exchanged after the flip, the parameters $\#UU$s and $\#DU$s have 
the same distribution on $\m_{n}$. In particular, $\v \m_{n}(UU) \v = 
\v \m_{n}(DU) \v$
(\htmladdnormallink{A004148}{http://www.research.att.com:80/cgi-bin/access.cgi/as/njas/sequences/eisA.cgi?Anum=A004148}
). 

Remark:\quad This bijection is clearly an involution on $\m_{n}$ and 
generalizes one on Dyck paths \cite{deutschInvol}. 
\qed 

In a full binary tree there is an obvious correspondence between 
non-first left leaves and right nodes: given such a leaf, travel 
(toward the root) to the first right node encountered.

\hspace*{40mm}\p[levelsep=30pt,nodesep=3pt]{\TR{\scr{$\circ$ }}}{ \p{\TR{\scr{$\circ$ }}}{ \TR{\scr{$\circ$ }} 
\p{\TR{\scr{$\circ$ }}}{\p{\TR{\scr{$\circ$ }}}{ \p{\TR{\scr{$\circ$ }}}{ \TR{\scr{$\circ$ }} 
\TR{\scr{$\circ$ }}} \TR{\scr{$\circ$ }}}\p{\TR{\scr{$\circ$ }}}{ \TR{\scr{$\circ$ }} \TR{\scr{$\circ$ }}}}
} \TR{\scr{$\circ$ }}}

% 1=step in x-direction, 2=step in y-direction, 3=upward diagonal step, 
% 4=downward diagonal step
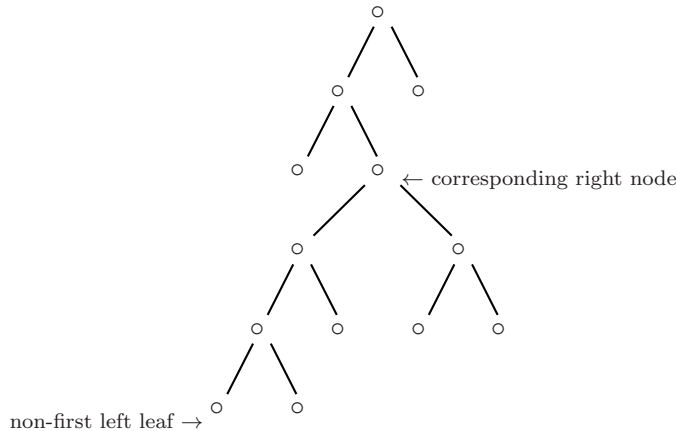
\vspace*{-35mm}

\Einheit=0.5cm
\[
\Label\o{\textrm{\scriptsize{$\leftarrow$ corresponding right 
node}}}(2,2.8)
\Label\o{\textrm{\scriptsize{non-first left leaf 
$\rightarrow$}}}(-9.5,-3.6)
\] 

\vspace*{2mm}

Note that this correspondence holds (vacuously) even for the trivial 
tree consisting of the root alone. 
By associating the root to the first left leaf, we can extend this 
correspondence to \{left leaves\} $\leftrightarrow $ \{right 
nodes\} $\cup$ \{root\}, except in the case of the trivial tree. This 
exception ultimately accounts for why many of our bijections need to 
increase the path length (tree weight) by 1.

Similarly, \{right leaves\} $\leftrightarrow $ \{left 
nodes\} $\cup$ \{root\} in all but the trivial tree. Applied to a Dyck 
path (all labels 0), these correspondences yield the obvious fact 
that \#peaks = \#valleys $ +1$ and the slightly less obvious fact 
that \#doublerises = \#doublefalls.

The next 3 bijections are all from $\m_{n}(UU)$
to $ \m_{n+1}(UD)$.
\begin{bij}
    $\phi: \m_{n}(UU) \rightarrow \m_{n+1}(UD)$.
\end{bij}
Recursive:
\begin{eqnarray*}
    \phi(\ve) & = & F  \\
    \phi(FR) & = & F\phi(R)  \\
    \phi(UF^{a}DR) & = & U\phi(R)DF^{a}\hspace*{30mm}a\ge 0  \\
    \phi(UF^{a}RF^{b}DS) & = &   U\phi(S)D \phi(F^{b}RF^{a-1}) 
    \hspace*{5mm} a\ge 1,\,b\ge 0;\ \textrm{$R$ starts $U$, ends $D$} 
\end{eqnarray*}

Explicit: Given a $UU$-free Motzkin $n$-path, its tree has positive 
labels on its left nodes. Flip the tree in the vertical and increment 
the root label by 1. Exchange labels on nonfirst left leaves and corresponding 
right nodes, and exchange the labels on the first 
leaf and the root. Now every left leaf has a positive label and 
the tree weight is incremented by 1. Take the corresponding path---a 
$UD$-free Motzkin $(n+1)$-path. The map is obviously reversible.

Example: 
% 1=step in x-direction, 2=step in y-direction, 3=upward diagonal step, 
% 4=downward diagonal step
\vspace*{2mm}
\Einheit=0.5cm
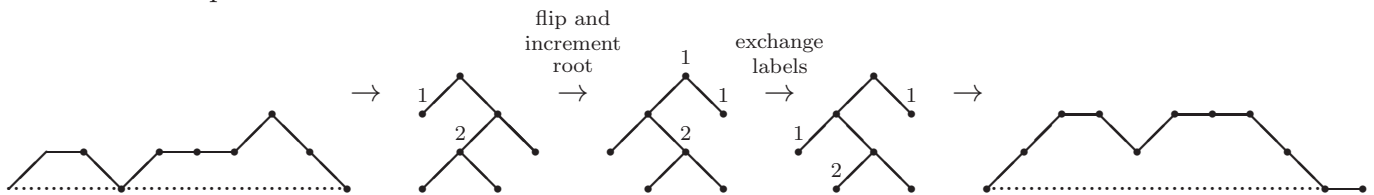
\[
\Pfad(-18,0),314311344\endPfad
\SPfad(-18,0),111111111\endSPfad
\Pfad(-7,2),344\endPfad
\Pfad(-7,0),334\endPfad
\Pfad(-6,1),4\endPfad
\Pfad(-2,1),34\endPfad
\Pfad(-1,2),34\endPfad
\Pfad(-1,0),34\endPfad
\Pfad(3,1),334\endPfad
\Pfad(4,2),4\endPfad
\Pfad(4,0),34\endPfad
\Pfad(8,0),3314311441\endPfad
\SPfad(8,0),111111111\endSPfad
\DuennPunkt(-16,1)
\DuennPunkt(-15,0)
\DuennPunkt(-14,1)
\DuennPunkt(-13,1)
\DuennPunkt(-12,1)
\DuennPunkt(-11,2)
\DuennPunkt(-10,1)
\DuennPunkt(-9,0)
\DuennPunkt(-7,2)
\DuennPunkt(-7,0)
\DuennPunkt(-6,1)
\DuennPunkt(-6,3)
\DuennPunkt(-5,2)
\DuennPunkt(-5,0)
\DuennPunkt(-4,1)
\DuennPunkt(-1,2)
\DuennPunkt(-2,1)
\DuennPunkt(-1,0)
\DuennPunkt(0,1)
\DuennPunkt(0,3)
\DuennPunkt(1,0)
\DuennPunkt(1,2)
\DuennPunkt(3,1)
\DuennPunkt(4,0)
\DuennPunkt(5,1)
\DuennPunkt(4,2)
\DuennPunkt(5,3)
\DuennPunkt(6,0)
\DuennPunkt(6,2)
\DuennPunkt(8,0)
\DuennPunkt(9,1)
\DuennPunkt(10,2)
\DuennPunkt(11,2)
\DuennPunkt(12,1)
\DuennPunkt(13,2)
\DuennPunkt(14,2) 
\DuennPunkt(15,2)
\DuennPunkt(16,1)
\DuennPunkt(17,0) 
\DuennPunkt(18,0)
\Label\o{\rightarrow}(-8.5,2)
\Label\o{\rightarrow}(-3.0,2)
\Label\o{\textrm{\scriptsize{flip and}}}(-3.0,4.0)
\Label\o{\textrm{\scriptsize{increment}}}(-3.0,3.4)
\Label\o{\textrm{\scriptsize{root}}}(-3.0,2.8)
\Label\o{\textrm{\scriptsize{exchange}}}(2.5,3.4)
\Label\o{\textrm{\scriptsize{labels}}}(2.5,2.8)
\Label\o{\rightarrow}(2.5,2)
\Label\o{\rightarrow}(7.5,2)
\Label\o{\textrm{{\scriptsize 1}}}(-7,2)
\Label\o{\textrm{{\scriptsize 2}}}(-6,1)
\Label\o{\textrm{{\scriptsize 2}}}(0,1)
\Label\o{\textrm{{\scriptsize 1}}}(0,3)
\Label\o{\textrm{{\scriptsize 1}}}(1,2)
\Label\o{\textrm{{\scriptsize 1}}}(3,1)
\Label\o{\textrm{{\scriptsize 2}}}(4,0)
\Label\o{\textrm{{\scriptsize 1}}}(6,2)
\] 

Consequence:\quad The parameters ``number of doublefalls'' ($DD$s) on $\m_{n}(UU)$ and 
``number of valleys'' ($DU$s) on 
$\m_{n+1}(UD)$ have 
the same distribution. This is because, after the tree flip and label 
exchange, a non-last right $0$-leaf ($DD$) becomes a right 0-node 
($DU$). In particular, $\v \m_{n}(UU,DD) \v = 
\v \m_{n+1}(UD,DU) \v$ or, more picturesquely, Motzkin $n$-paths 
containing no long slanted segments are equinumerous with 
Motzkin $(n+1)$-paths containing no sharp turns
(\htmladdnormallink{A004149}{http://www.research.att.com:80/cgi-bin/access.cgi/as/njas/sequences/eisA.cgi?Anum=A004149}). \qed 
%\newpage 

\begin{bij}
     $\phi: \m_{n}(UU) \rightarrow \m_{n+1}(UD)$.
\end{bij}
Recursive:
\begin{eqnarray*}
    \phi(\ve) & = & F  \\
    \phi(FR) & = & F\phi(R)  \\
    \phi(UDR) & = & U\phi(R)D  \\
    \phi(UFRDS) & = &   U\phi(R)D \phi(S) 
\end{eqnarray*}

Explicit: Given a $UU$-free Motzkin $n$-path, its tree has positive 
labels on its left nodes. Consider the labels as counting tokens 
(flatsteps) stored at their location. \\
Step 1.\ Add a token to the root. \\
Step 2.\ Transfer one token from each left node and from the root to 
its corresponding right leaf. (Except do nothing if the tree consists 
of the root alone.) \\
Step 3.\ For each left $0$-leaf, transfer the subtree of its sibling 
vertex (including the label on the sibling vertex) to this left leaf. 
Note that after the transfer, the left $0$-leaf may become a node or 
a positive leaf but will no longer be a 0-leaf. Also, the sibling in 
question becomes a right 0-leaf and all other right leaves are 
positive due to Step 2. The weight has been increased by 1 and all 
left leaves are now positive, so the resulting path is indeed in 
$\m_{n+1}(UD)$. The map is reversible: the original left 0-leaves are 
recovered as the siblings of right 0-leaves in the image.

Example: (same $UU$-free path as in the previous example)

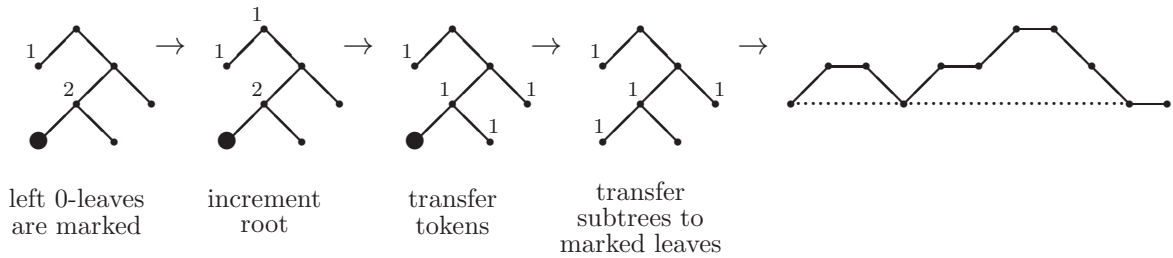
\Einheit=0.5cm
\[
\Pfad(-15,0),334\endPfad
\Pfad(-15,2),34\endPfad
\Pfad(-14,1),4\endPfad
\DuennPunkt(-15,2)
\DickPunkt(-15,0)
\DuennPunkt(-14,1)
\DuennPunkt(-14,3)
\DuennPunkt(-13,0)
\DuennPunkt(-13,2)
\DuennPunkt(-12,1)
\Pfad(-10,0),334\endPfad
\Pfad(-10,2),34\endPfad
\Pfad(-9,1),4\endPfad
\DuennPunkt(-10,2)
\DickPunkt(-10,0)
\DuennPunkt(-9,1)
\DuennPunkt(-9,3)
\DuennPunkt(-8,0)
\DuennPunkt(-8,2)
\DuennPunkt(-7,1)
\Pfad(-5,0),334\endPfad
\Pfad(-5,2),34\endPfad
\Pfad(-4,1),4\endPfad
\DuennPunkt(-5,2)
\DickPunkt(-5,0)
\DuennPunkt(-4,1)
\DuennPunkt(-4,3)
\DuennPunkt(-3,0)
\DuennPunkt(-3,2)
\DuennPunkt(-2,1)
\Pfad(0,0),334\endPfad
\Pfad(0,2),34\endPfad
\Pfad(1,1),4\endPfad
\DuennPunkt(0,2)
\DuennPunkt(0,0)
\DuennPunkt(1,1)
\DuennPunkt(1,3)
\DuennPunkt(2,0)
\DuennPunkt(2,2)
\DuennPunkt(3,1)
\Pfad(5,1),3143131441\endPfad
\SPfad(5,1),111111111\endSPfad
\DuennPunkt(5,1)
\DuennPunkt(6,2)
\DuennPunkt(7,2)
\DuennPunkt(8,1)
\DuennPunkt(9,2)
\DuennPunkt(10,2)
\DuennPunkt(11,3)
\DuennPunkt(12,3)
\DuennPunkt(13,2)
\DuennPunkt(14,1) 
\DuennPunkt(15,1) 
\Label\o{\rightarrow}(-11.5,2)
\Label\o{\rightarrow}(-6.5,2)
\Label\o{\rightarrow}(-1.5,2)
\Label\o{\rightarrow}(4,2)
\Label\o{\textrm{\footnotesize{left 0-leaves}}}(-14.0,-2.0)
\Label\o{\textrm{\footnotesize{are marked}}}(-14.0,-2.7)
\Label\o{\textrm{{\scriptsize 1}}}(-15.2,1.9)
\Label\o{\textrm{{\scriptsize 2}}}(-14.2,0.9)
\Label\o{\textrm{\footnotesize{increment}}}(-9.0,-2.0)
\Label\o{\textrm{\footnotesize{root}}}(-9.0,-2.7)
\Label\o{\textrm{{\scriptsize 1}}}(-9.2,2.9)
\Label\o{\textrm{{\scriptsize 1}}}(-10.2,1.9)
\Label\o{\textrm{{\scriptsize 2}}}(-9.2,0.9)
\Label\o{\textrm{\footnotesize{transfer}}}(-4.0,-2.0)
\Label\o{\textrm{\footnotesize{tokens}}}(-4.0,-2.7)
\Label\o{\textrm{{\scriptsize 1}}}(-5.2,1.9)
\Label\o{\textrm{{\scriptsize 1}}}(-4.2,0.9)
\Label\o{\textrm{{\scriptsize 1}}}(-2.9,-0.1)
\Label\o{\textrm{{\scriptsize 1}}}(-1.9,0.9)
\Label\o{\textrm{\footnotesize{transfer}}}(1.0,-1.8)
\Label\o{\textrm{\footnotesize{subtrees to}}}(1.0,-2.5)
\Label\o{\textrm{\footnotesize{marked leaves}}}(1.0,-3.2)
\Label\o{\textrm{{\scriptsize 1}}}(-0.2,1.9)
\Label\o{\textrm{{\scriptsize 1}}}(0.8,0.9)
\Label\o{\textrm{{\scriptsize 1}}}(-0.1,-0.1)
\Label\o{\textrm{{\scriptsize 1}}}(3.1,0.9)
\] 

Consequence:\quad Define a low peak in a Motzkin path to be a peak 
whose downstep returns the path to ground level, and the final descent 
to be the one that terminates the path (assumed empty if the path ends 
with a flatstep). Then the parameters ``\#low peaks'' on $\m_{n}(UU)$ 
and ``length of final descent'' on $\m_{n+1}(UD)$ have the same distribution. 
In particular, $\v \{P \in \m_{n}(UU): P \textrm{ has no low peaks}\} 
\v = \v \m_{n}(UD) \v$. This is because the image ends $F$ (i.e. has 
final descent 0) $\Leftrightarrow$ the original path has no low 
peaks. Deleting this final $F$ is a bijection to $\m_{n}(UD)$.  \qed 

\begin{bij}
     $\phi: \m_{n}(UU) \rightarrow \m_{n+1}(UD)$.
\end{bij}
Recursive:
\begin{eqnarray*}
    \phi(\ve) & = & F  \\
    \phi(FR) & = & F\phi(R)  \\
    \phi(UDR) & = & U\phi(R)D  \\
    \phi(UFRDS) & = &   U\phi(S)D \phi(R) 
\end{eqnarray*}

Explicit: Given a $UU$-free Motzkin $n$-path, pass to its corresponding 
tree, flip tree in vertical and add a token to the root. Now the root 
and right nodes all have positive labels. Transfer one token from each to the 
corresponding left leaf (do nothing if there are no edges). Here 
again, all left leaves are now positive, giving a $UD$-free Motzkin 
$(n+1)$-path.

Example: (same $UU$-free path as in the previous two examples)

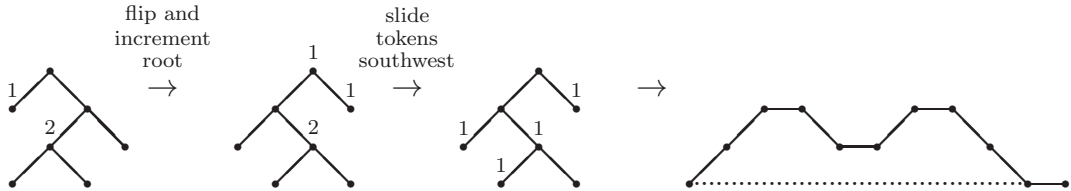
\hspace*{-50mm}

\Einheit=0.5cm
\[
\Pfad(-14,2),344\endPfad
\Pfad(-14,0),334\endPfad
\Pfad(-13,1),4\endPfad
\DuennPunkt(-14,2)
\DuennPunkt(-14,0)
\DuennPunkt(-13,1)
\DuennPunkt(-13,3)
\DuennPunkt(-12,2)
\DuennPunkt(-12,0)
\DuennPunkt(-11,1)
\Pfad(-8,1),34\endPfad
\Pfad(-7,2),34\endPfad
\Pfad(-7,0),34\endPfad
\DuennPunkt(-7,2)
\DuennPunkt(-8,1)
\DuennPunkt(-7,0)
\DuennPunkt(-6,1)
\DuennPunkt(-6,3)
\DuennPunkt(-5,0)
\DuennPunkt(-5,2)
\Pfad(-2,1),334\endPfad
\Pfad(-1,2),4\endPfad
\Pfad(-1,0),34\endPfad
\DuennPunkt(-2,1)
\DuennPunkt(-1,0)
\DuennPunkt(0,1)
\DuennPunkt(-1,2)
\DuennPunkt(0,3)
\DuennPunkt(1,0)
\DuennPunkt(1,2)
\Pfad(4,0),3314131441\endPfad
\SPfad(4,0),111111111\endSPfad
\DuennPunkt(4,0)
\DuennPunkt(5,1)
\DuennPunkt(6,2)
\DuennPunkt(7,2)
\DuennPunkt(8,1)
\DuennPunkt(9,1)
\DuennPunkt(10,2) 
\DuennPunkt(11,2)
\DuennPunkt(12,1)
\DuennPunkt(13,0) 
\DuennPunkt(14,0)
\Label\o{\rightarrow}(-10.0,2)
\Label\o{\textrm{\scriptsize{flip and}}}(-10.0,4.0)
\Label\o{\textrm{\scriptsize{increment}}}(-10.0,3.4)
\Label\o{\textrm{\scriptsize{root}}}(-10.0,2.8)
\Label\o{\textrm{\scriptsize{slide}}}(-3.5,4.0)
\Label\o{\textrm{\scriptsize{tokens}}}(-3.5,3.4)
\Label\o{\textrm{\scriptsize{southwest}}}(-3.5,2.8)
\Label\o{\rightarrow}(-3.5,2)
\Label\o{\rightarrow}(3.0,2)
\Label\o{\textrm{{\scriptsize 1}}}(-14,2)
\Label\o{\textrm{{\scriptsize 2}}}(-13,1)
\Label\o{\textrm{{\scriptsize 2}}}(-6,1)
\Label\o{\textrm{{\scriptsize 1}}}(-6,3)
\Label\o{\textrm{{\scriptsize 1}}}(-5,2)
\Label\o{\textrm{{\scriptsize 1}}}(-2,1)
\Label\o{\textrm{{\scriptsize 1}}}(-1,0)
\Label\o{\textrm{{\scriptsize 1}}}(0,1)
\Label\o{\textrm{{\scriptsize 1}}}(1,2)
\] 
Consequence:\quad The parameters $\#UFU$s on $\m_{n}(UU)$ and $\#DU$s on 
$\m_{n+1}(UD)$ have 
the same distribution. This is because a $UFU$ occurs for each 
left node with label 1. After the flip, it becomes a right node with label 1.
The required transfer of tokens makes it a right 0-node, and a $DU$ in 
the image path. \qed 

Next, we consider a bijection that increases plateau lengths in 
valley-free paths. Recall that a $U\underbrace{F\ldots F}_{k \ge 1}D$ 
sequence is a plateau of length $k$ (as is 
$\underbrace{F\ldots F}_{k \ge 0}$ 
if it is the entire path). Let $MPL$ 
denote minimum plateau length 
(\htmladdnormallink{A064645}{http://www.research.att.com:80/cgi-bin/access.cgi/as/njas/sequences/eisA.cgi?Anum=A064645})
in a path (taken to be 0 if there are no plateaus).

\begin{bij}
     $\phi: \m_{n}(DU) \rightarrow \m_{n+1}(UD)$.
\end{bij}
Recursive:
\begin{eqnarray*}
    \phi(\ve) & = & F  \\
    \phi(FR) & = & F\phi(R)  \\
    \phi(URD) & = & U\phi(R)D  \\
    \phi(URDFS) & = &   U\phi(R)D \phi(S) 
\end{eqnarray*}

Explicit: A $DU$-free path gives a tree in which all right nodes have 
a positive label. Add a token to the root and (unless it's the 
root-only tree) transfer 
a token from the root and from each right node to its corresponding left 
leaf. Now all left leaves have a positive label, giving a $UD$-free 
Motzkin $(n+1)$-path. The map is clearly reversible.

Example: 

\hspace*{100mm} \p[levelsep=30pt,nodesep=3pt]{\TR{1}} 
                               {\TR{\scr{1}}           \p{\TR{1}}
                                           {\p{\TR{$\circ$}}            
			      {\TR{\scr{2}}           \TR{$\circ$}} \p{\TR{\scr{2}}} 
			      {\TR{$\circ$} \TR{$\circ$}}}}

\vspace*{-20mm}
\Einheit=0.5cm
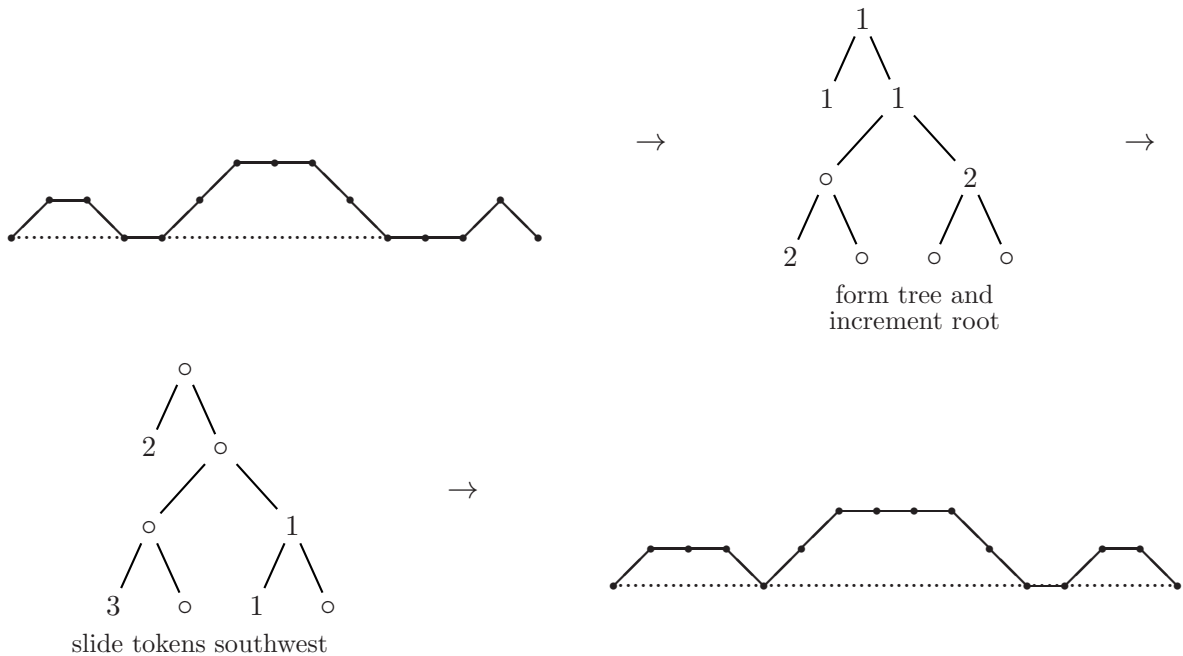
\[
\Pfad(-15,0),31413311441134\endPfad
\SPfad(-15,0),111\endSPfad
\SPfad(-11,0),111111\endSPfad
\DuennPunkt(-15,0)
\DuennPunkt(-14,1)
\DuennPunkt(-13,1)
\DuennPunkt(-12,0)
\DuennPunkt(-11,0)
\DuennPunkt(-10,1)
\DuennPunkt(-9,2)
\DuennPunkt(-8,2)
\DuennPunkt(-7,2)
\DuennPunkt(-6,1)
\DuennPunkt(-5,0)
\DuennPunkt(-4,0)
\DuennPunkt(-3,0)
\DuennPunkt(-2,1)
\DuennPunkt(-1,0)
\Label\o{\rightarrow}(2,2)
\Label\o{\rightarrow}(15,2)
\Label\o{\textrm{{\footnotesize form tree and}}}(9,-2)
\Label\o{\textrm{{\footnotesize increment root}}}(9,-2.7)
\]

\hspace*{10mm} \p[levelsep=30pt,nodesep=3pt]{\TR{$\circ$}} {\TR{\scr{2}} \p{\TR{$\circ$}} 
{\p{\TR{$\circ$}} {\TR{\scr{3}} \TR{$\circ$}} \p{\TR{\scr{1}}} {\TR{\scr{1}} \TR{$\circ$}}}}

\vspace*{-20mm}
\Einheit=0.5cm
\[
\Pfad(1,0),311433111441314\endPfad
\SPfad(1,0),111111111111\endSPfad
\SPfad(13,0),111\endSPfad
\DuennPunkt(1,0)
\DuennPunkt(2,1)
\DuennPunkt(3,1)
\DuennPunkt(4,1)
\DuennPunkt(5,0)
\DuennPunkt(6,1)
\DuennPunkt(7,2)
\DuennPunkt(8,2)
\DuennPunkt(9,2)
\DuennPunkt(10,2)
\DuennPunkt(11,1)
\DuennPunkt(12,0) 
\DuennPunkt(13,0)
\DuennPunkt(14,1)
\DuennPunkt(15,1) 
\DuennPunkt(16,0)
\Label\o{\rightarrow}(-3,2)
\Label\o{\textrm{\footnotesize{slide tokens southwest}}}(-10,-2)
\]

Consequence:\quad The parameters $MPL+1$ on $\m_{n}(DU)$ and $MPL$ on 
$\m_{n+1}(UD) $ have the same distribution. This follows since each 
peak (left 0-leaf) becomes a plateau of length 1 and each existing 
plateau has its length increased by 1. \qed

Bijection 4 followed by the inverse of Bijection 2, that is, 
$\phi_{2}^{-1}\phi_{4}$, shows that the parameters $\#UFU$s and 
$\#DD$s are equidistributed on $\m_{n}(UU)$. In fact, this bijection 
is an involution on $\m_{n}(UU)$ that interchanges occurrences of 
$UFU$ and $DD$. We conclude with a simple explicit description of 
$\phi_{2}^{-1}\phi_{4}$. For this purpose say an upstep is critical 
if it is followed by an $F$. A strict Motzkin path is one that 
starts $U$ and ends $D$. Given a $UU$-free path, if the path 
segment strictly between a critical $U$ and its matching $D$ is level 
(all $F$s), leave it alone. Otherwise, it has the form $F^{a}SF^{b}$ 
with $a\ge 1,\,b \ge 0,\,S$ strict. Replace 
it by $F^{b}SF^{a-1}$. The result is independent of the 
order in which critical $U$s are processed and is again $UU$-free. 
This map 
is clearly an involution and, because of the restriction to $UU$-free 
paths, $DD$s in the image path correspond one-to-one to $UFU$s in the 
original.

\begin{center}
{\normalsize \textbf{Acknowledgement}}
\end{center}
\vspace*{-5mm}
I thank Emeric Deutsch for several helpful comments on an early 
draft of this paper.

\end{document}